\documentclass[manuscript, screen]{acmart}
\pdfoutput=1
\usepackage[mode=multiuser,draft]{fixme} 
\FXRegisterAuthor{cc}{envcc}{CC}
\FXRegisterAuthor{gk}{envgk}{GK}
\fxusetheme{color}
\AtBeginDocument{%
  }

\setcopyright{acmlicensed}
\copyrightyear{2024}
\acmYear{2024}
\acmDOI{XXXXXXX.XXXXXXX}

\usepackage{tikz,tikzscale}
\usepackage{graphicx}

\usetikzlibrary{calc}
\usetikzlibrary{shadings}
\usetikzlibrary{shapes.geometric}
\usetikzlibrary{shadows.blur}
\usetikzlibrary{positioning}
\usetikzlibrary{backgrounds}
\usetikzlibrary{matrix}

\usepackage{pgfplots}
\usepgfplotslibrary{groupplots}

\usepackage{algorithm}
\usepackage{algpseudocodex}

\usepackage{xspace}
\usepackage{bold-extra}
\usepackage[most]{tcolorbox}

\usepackage{amsmath}
\usepackage{amsfonts}

\usepackage{stmaryrd}
\SetSymbolFont{stmry}{bold}{U}{stmry}{m}{n}

\usepackage{braket,amsfonts}

\usepackage{array}

\usepackage{mathtools}

\usepackage{arydshln}

\usepackage{multirow}

\begin{document}

\title{Multilevel Interior Penalty Methods on GPUs}

\author{Cu Cui}
\authornote{Submitted to the editors April 03, 2024.}
\email{cu.cui@iwr.uni-heidelberg.de}
\author{Guido Kanschat}
\email{kanschat@uni-heidelberg.de}
\affiliation{%
  \department[1]{Interdisciplinary Center for Scientific Computing (IWR)}
  \institution{Heidelberg University}
  \streetaddress{Im Neuenheimer Feld 205}
  \city{Heidelberg}
  \country{Germany}
  \postcode{69120}
}

\begin{abstract}
  We present a matrix-free multigrid method for high-order discontinuous Galerkin (DG) finite element methods with GPU acceleration. A performance analysis is conducted, comparing various data and compute layouts. Smoother implementations are optimized through localization and fast diagonalization techniques. Leveraging conflict-free access patterns in shared memory, arithmetic throughput of up to 39\% of the peak performance on Nvidia A100 GPUs are achieved. Experimental results affirm the effectiveness of mixed-precision approaches and MPI parallelization in accelerating algorithms. Furthermore, an assessment of solver efficiency and robustness is provided across both two and three dimensions, with applications to Poisson problems.
\end{abstract}




\maketitle

\section{Introduction}
\label{sec:intro}

Over the past three decades, Discontinuous Galerkin (DG) methods have garnered considerable attention in various fields.
These inherently high-order methods are favored for achieving superior accuracy per degree of freedom compared to their lower-order counterparts. Furthermore, their high arithmetic intensity makes them ideally suited for deployment on graphics processing units (GPUs) and the GPU-accelerated architectures prevalent in current and upcoming supercomputers~\cite{klockner2009nodal,vermeire2017utility,brown2018ceed}.

However, traditional matrix-based approaches often become impractical for high polynomial degrees due to escalating computational costs and memory demands. A more viable solution is to replace matrix-vector products with on-the-fly evaluations of discretized differential operators in a matrix-free manner. Utilizing sum-factorization techniques on tensor-product elements~\cite{orszag1979spectral,patera1984spectral}, these matrix-free operators can be executed with $O(k^{d+1})$ operations and require $O(k^d)$ memory~\cite{pazner2018approximate}. Moreover, recent work has demonstrated near-peak performance on GPUs~\cite{brown2018ceed,swirydowicz2019acceleration}, positioning them as an attractive choice for efficient computational implementations.

Recent studies have focused on matrix-free implementations tailored to GPUs~\cite{modave2016gpu,remacle2016gpu,ljungkvist2017matrix}, with particular emphasis on matrix-free multigrid methods employing point Jacobi and Chebyshev smoothing~\cite{KronbichlerLjunqkvist19}, and matrix-free tensor-product approximations for block preconditioners in discontinuous Galerkin discretizations~\cite{franco2020high}. Yet, these methods often face performance limitations due to the reliance on point or simple block smoothers. Given the role of a robust smoother to achieve fast multigrid convergence for high-order discretizations~\cite{kanschat2008robust,kanschat2015multigrid,Cui2023}, there is a need to develop efficient GPU algorithms that transcend these constraints.

In this paper, we develop a matrix-free multigrid preconditioner with a vertex-patch smoother on GPUs. A patch-wise integration approach for matrix-free evaluation of finite element operators with constant coefficients on a Cartesian mesh is considered. This approach not only facilitates structured data access but also reduces the overall arithmetic and data access operations. The resulting preconditioner exhibits robustness across varying mesh sizes and polynomial degrees. Notably, for higher-order elements, only few steps are required to achieve convergence, rendering it nearly equivalent to a direct solver. Furthermore, when combined with fast diagonalization techniques, the computational effort for overlapping vertex-patch smoothing becomes comparable to that of a matrix-free operator application.

This paper is organized as follows: In Section~\ref{eq:model}, we briefly describe the model problem and its discretization, with a specific focus on the implementation of the matrix-free operators. Matrix-free smoothers and preconditioners are developed in Section~\ref{sec:GMG}. Section~\ref{sec:Imp_s} details the GPU implementation, and Section~\ref{sec:opt} delves into the algorithmic optimizations. The paper then explores MPI parallelization in Section~\ref{sec:Imp_m} and Section~\ref{sec:con} concludes this work.

\section{Interior Penalty Methods}
\label{sec:model}

In this study, we consider the approximation of solutions to the Poisson equation
\begin{gather}\label{eq:model}
\begin{aligned}
    -\Delta u &= f & \text{ in } & \Omega \\
    u &{} = 0 & \text{ on } &\partial\Omega
\end{aligned}
\end{gather}
subject to Dirichlet boundary conditions on a domain $\Omega \subset \mathbb{R}^{d}$ with dimension $d=2,3$. 
The approach is shown for the Laplacian, but it can be applied to other symmetric, separable differential operators. 

\subsection{Discretization}
The finite element space is  chosen such that the restriction of $V_h$ to a single cell $K \in \mathcal{T}_h$ is $\mathbb{Q}_k$, the space of mapped tensor product polynomials of degree $k$.
To define bilinear forms associated with DG discretizations of the Laplacian, we introduce the following notation for jumps and mean values of discontinuous functions on interior faces:
\begin{equation}\label{eq:jump_mean}
    \{u\} = \frac{u^+ + u^-}{2}, \quad \llbracket u \rrbracket = (u^+ + u^-)\boldsymbol{n}^+,
\end{equation}
and on boundary faces:
\begin{equation}\label{eq:jump_boundary}
    \{u\} = u, \quad \llbracket u \rrbracket = u.
\end{equation}
The set of faces between mesh cells and at the boundary is denoted as $\mathcal{F}_h$. We adopt the symmetric interior penalty method (as detailed in~\cite{arnold1982interior}) and introduce the associated bilinear form as

\begin{equation}
\label{eq:bilinear_form}
\begin{aligned}
    a_h(u_h,v_h) &=  \sum_{K \in \mathcal{T}_h} \int_{K} \nabla u_h \cdot \nabla v_h \,d\boldsymbol{x} \\ &
     +\sum_{F \in \mathcal{F}_h} \int_{F}\left(\gamma_e \llbracket u_h \rrbracket \cdot \llbracket v_h \rrbracket-\{\nabla u_h\} \cdot \llbracket v_h \rrbracket-\llbracket u_h \rrbracket \cdot\{\nabla v_h\}\right) \,d \sigma(\boldsymbol{x}),
\end{aligned}
\end{equation}
where the penalty parameter $\gamma_e$ takes the form of~\cite{kanschat2003discontinuous} 
\begin{equation}
    \gamma_e = k(k+1)\left( \frac{1}{h^+} + \frac{1}{h^-} \right).
\end{equation}
Finally, the interior penalty discretization of the model problem~\eqref{eq:model} reads: find $u_h \in V_h$ such that
\begin{equation}\label{eq:weak_form}
     a_h(u_h,v_h) = \sum_{K \in \mathcal{T}_h} \int_{K} f v_h \,d\boldsymbol{x} 
     \quad \forall v_h \in V_h,
\end{equation}
The goal of the preconditioner outlined here is to solve~\eqref{eq:weak_form} efficiently. To this end, we refer to it in matrix form
\begin{equation}\label{eq:mat_form}
    Ax = b.
\end{equation}

\subsection{Matrix-free implementation}\label{sec:mf_op}

Constructing the matrix defined in~\eqref{eq:mat_form} through explicit assembly typically costs $\mathcal{O}(k^{3d})$ operations when employing a straightforward algorithm. Additionally, the storage requirements for these matrices amount to $\mathcal{O}(k^{2d})$ memory. These computational and memory demands can become prohibitively expensive, particularly in simulations involving high polynomial degrees. This challenge becomes more pronounced on GPUs, where the system matrices may exceed the device memory capacity. To address these issues, we choose a matrix-free implementation, replacing matrix-vector products with on-the-fly integral evaluations~\cite{KronbichlerKormann19,KronbichlerLjunqkvist19,swirydowicz2019acceleration}.

The matrix-free evaluation of the finite element operator $Ax$ in~\eqref{eq:mat_form} is achieved by looping over all cells and faces appearing in the operator~\eqref{eq:bilinear_form}. Interpolating the solution from the coefficient values to quadrature points and performing the summation for quadrature requires $\mathcal{O}(k^{2d})$ arithmetic operations. By exploiting the tensor product form of basis functions and quadrature formula, we can reduce the complexity to $\mathcal{O}(dk^{d+1})$ by means of sum factorization~\cite{swirydowicz2019acceleration,vos2010h}.

For the case of constant coefficients and Cartesian geometry, the cell matrix can be represented as a sum of Kronecker products, i.e. 
\begin{equation*}
    A_j = L_1 \otimes M_0 + M_1 \otimes L_0
\end{equation*}
in 2D and 
\begin{equation*}
 A_j = L_2 \otimes M_1 \otimes M_0 + M_2 \otimes L_1 \otimes M_0 + M_2 \otimes M_1 \otimes L_0
\end{equation*}
in 3D, respectively, where $L_d$ and $M_d$ are one dimensional stiffness matrices and mass matrices, respectively. More details can be found in~\cite{witte2021fast}. This specific optimization involves fewer sum-factorization operations than numerical integration.

In the context of face terms, a face-wise integration approach was explored in the work~\cite{KronbichlerKormann19}. This method offers the advantage of simultaneously computing the flux term for both sides of a face, leading to a reduction in the overall number of arithmetic and data access operations at quadrature points. Conversely, the cell-wise evaluation of face integrals computes all contributions to a cell in a single pass, resulting in a more structured access to vectors organized on a cell-wise basis.

\begin{figure}[tp]
\centering
\includegraphics[width=.65\textwidth]{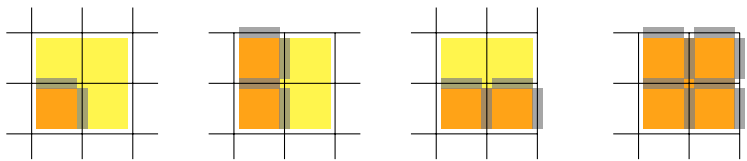}
\caption{Compute pattern for patch-wise integrals in 2D. Orange indicates read/write access by the current patch, yellow read only. Gary indicates a contribution of marked faces to orange cell(s). From left to right: center patch, top boundary patch, right boundary patch and corner patch.}
\label{fig:patch_int}
\Description[Fully described in the text.]{}
\end{figure}

In this work, we present a patch-wise integration pattern that merges these methodologies for enhanced computational efficiency. This approach allows simultaneous computation of integrals for partial cells and faces within a patch, demonstrated in Figure~\ref{fig:patch_int}. Based on our tests, the patch-wise pattern demonstrates better performance owing to its notable advantages in terms of accessing contiguous memory addresses, a critical aspect for achieving coalesced memory access. Moreover, the benefits of the patch-wise operation extend to both the smoothing operator in the V-cycle and the grid transfer operator, as detailed later in this paper. The consistent use of a unified data structure for these operations not only results in memory savings but also streamlines data management. It is essential to underscore that our optimization is specifically tailored for axis-aligned meshes. The primary objective of this work is the efficient implementation of mathematically efficient algorithms on GPUs, aiming to achieve optimal performance.

\section{Multigrid Method}\label{sec:GMG}

To establish the multilevel preconditioner, as introduced in~\cite{gopalakrishnan2003multilevel}, we assume that we have a hierarchy of meshes
\begin{gather*}
    \mathbb{T}_0 \sqsubset \mathbb{T}_1 \sqsubset\dots\sqsubset\mathbb{T}_L,
\end{gather*}
where each cell within $\mathbb{T}_\ell$ with $\ell \geq 1$ results from the regular refinement of a cell in $\mathbb{T}_{\ell-1}$ , leading to the generation of $2^d$ children.
Throughout our experiments, we choose $\mathbb{T}_0$ as the refinement of a single cell into $2^d$ cells.
Associated with these subdivisions, we have a hierarchy of nested finite element spaces
\begin{gather*}
    V_0 \subset V_1\subset \dots\subset V_L,
\end{gather*}
Between the spaces $V_\ell$, we define the prolongation operator $I_{\ell}^{\uparrow}$ from $V_{\ell}$ to $V_{\ell+1}$ as the canonical embedding and the restriction operator $I_{\ell}^{\downarrow}$ as the $L_2$-projection from $V_{\ell+1}$ to $V_{\ell}$. These are standard choices for multilevel methods in finite element computations.


We will consider the V-cycle as the most basic multigrid cycle. The multigrid preconditioner $P_\ell^{-1}$ is defined recursively as follows: on the coarsest level we let $P_0^{-1} = A_0^{-1}$. For level $\ell \geq 1$, define the action of $P^{-1}_\ell$ on a vector $b_\ell$ by
\begin{itemize}
    \item[(1)] Pre-smoothing: 
    \begin{equation*}
        x_\ell \gets S_\ell(x_\ell,b_\ell)
    \end{equation*} 
    \item[(2)] Coarse grid correction: 
    \begin{equation*}
        x_\ell \gets x_\ell + I_{\ell-1}^\uparrow P^{-1}_{\ell-1} I_{\ell-1}^\downarrow (b_\ell - A_\ell x_\ell)
    \end{equation*}
    \item[(3)] Post-smoothing:
    \begin{equation*}
        x_\ell \gets S_\ell(x_\ell,b_\ell)
    \end{equation*}
    \item[] and let $P^{-1}_{\ell-1}b_\ell = x_\ell$.
\end{itemize}
The operators $S_\ell$ represent the smoothers on level $\ell$, and their detailed description is provided in the subsequent subsection.
For the sake of simplicity in our discussion, we consistently assume a single pre- and post-smoothing step.

\subsection{Vertex-patch Smoothers}

We use a multiplicative vertex-patch smoother,
which is an overlapping domain decomposition method.
Each subdomain consists of all cells sharing a single vertex.
These smoothers were originally designed for spaces such as $H^{\text{div}}$ and $H^{\text{curl}}$ in~\cite{ArnoldFalkWinther97Hdiv,ArnoldFalkWinther00}.
They have been successfully applied to address problems involving Stokes~\cite{kanschat2015multigrid} and Darcy–Stokes–Brinkman~\cite{kanschat2017geometric} problems.

\begin{figure}[tp]
\centering
\includegraphics[width=.65\textwidth]{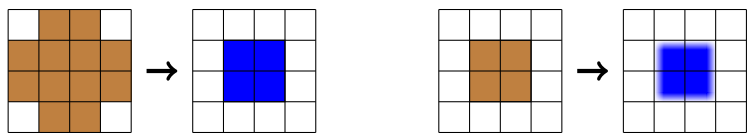}
\caption{Left: domain of dependence $\overline V_j$ (brown) and the support of $V_j$ (blue) of the vertex patch solver for the discontinuous Galerkin method. Right: the same for the continuous method where fading indicates that degrees on the boundary of the patch are not part of $V_j$.}
\label{fig:smoother}
\Description[Fully described in the text.]{}
\end{figure}

Let $V_j$ be the subspace of finite element functions with support on the patch around a vertex.
The local correction of a current approximate solution $x$ then consists of two steps, first computing the residual
\begin{gather}
    \label{eq:residual}
    r \gets b-Ax,
\end{gather}
and then the local solver
\begin{gather}
\label{eq:local-solver}
    x \gets x+R_j^T A_j^{-1}R_j r.
\end{gather}
Since memory access is usually a limiting factor on modern hardware, the computation of the residual should not be separated from the local solvers.
The restriction operator $R_j$ selectively extracts coefficients from the vector $r$ which are associated with basis functions in $V_j$. Similarly, the matrix $A_j$ is obtained by restricting the definition of the matrix in~\eqref{eq:mat_form} to the basis functions spanning $V_j$.
This means that only those entries of the vector $x$ in~\eqref{eq:residual} are required for the local solver of patch $j$ which couple through the matrix $A$ with the subspace $V_j$.
We refer to this subspace as well as its support as the domain of dependence $\overline V_j$.
This is illustrated in Figure~\ref{fig:smoother} on the left, where the support of $V_j$ is indicated in blue and the domain of dependence in brown.
We refer to this local solver as the \textit{full} kernel.
A challenge with this approach is that we lose the tensor product structure, leading to an irregular data structure.
Moreover, data has to be read from twice as many cells as contained in the vertex patch itself, increasing the load on the memory bus considerably.
Finally, data on the cells adjacent to a patch cannot be updated in parallel, limiting opportunities for parallelism.

We compare to the situation with continuous finite elements in Figure~\ref{fig:smoother} on the right. There, it is possible to solve on each vertex patch with homogeneous Dirichlet boundary conditions, such that $V_j$ is the subspace containing all ``interior'' degrees of freedom of the patch, while $\overline V_j$ contains all degrees of freedom of the same patch and not more. Hence, less and better structured data is read to compute local residuals.

We mimic this effect for discontinuous finite elements by using basis functions based on Lagrange interpolation including points on the boundaries of the mesh cells as in Figure~\ref{fig:dof_layout} on the left. Omitting the functions associated to boundary interpolation points from $V_j$, the functions in this space will vanish on the boundary. For the space $\overline V_j$, we use only the functions on patch $j$, such that the domain of dependence is the same as for continuous elements. We refer to this local solver as \textit{Dirichlet} kernel.
\begin{figure}[tp]
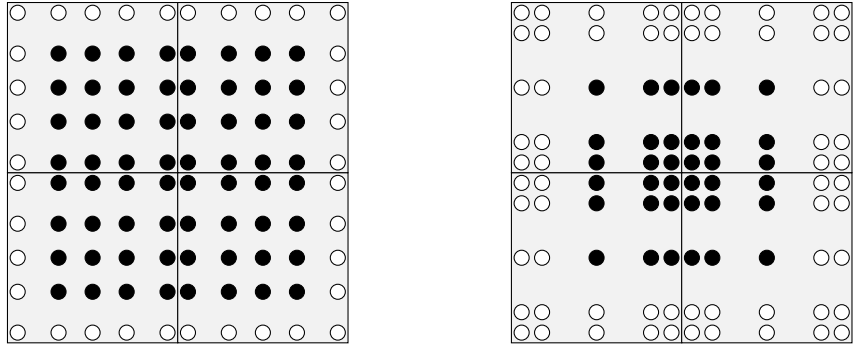

\centering
\includegraphics[width=.3\textwidth]{Figures/lagrange.tikz}
\hspace{2cm}
\includegraphics[width=.3\textwidth]{Figures/hermite.tikz}
\caption{Data access pattern for the Dirichlet and clamped kernels, respectively, from left to right with polynomial degree $k = 4$. Disks indicate read/write access corresponding to the space $V_j$, circles read only, corresponding to $\overline V_j/V_j$.}
\label{fig:hermite_type}
\Description[Fully described in the text.]{}
\end{figure}

We point out though, that doing so, the patch solver is not using the correct residual anymore, since the last face term in the bilinear form~\eqref{eq:bilinear_form} does not vanish in this case.
In order to address this inconsistency, we introduce a new local solver. To this end, we use a local basis which involves Hermite interpolation at cell boundaries as in Figure~\ref{fig:hermite_type} on the right.
Hence, we can easily restrict the space $V_j$ to functions $v$ with $v=0$ and $\partial_n v=0$ on the boundary of the patch.
All boundary terms in the bilinear form~\eqref{eq:bilinear_form} vanish for such test functions, such that the residual does not couple to neighboring cells.
Consequently, this approach leads to the exact residual on the space $V_j$ using the domain of dependency of the continuous method.
We refer to this local solver as \textit{clamped} kernel in reference to the same type of boundary conditions for the plate equation.
The data access pattern for a Hermite-type basis with polynomial degree $k=4$ is illustrated in Figure~\ref{fig:hermite_type}.

\begin{figure}[tp]
\centering
    \includegraphics[width=0.85\textwidth]{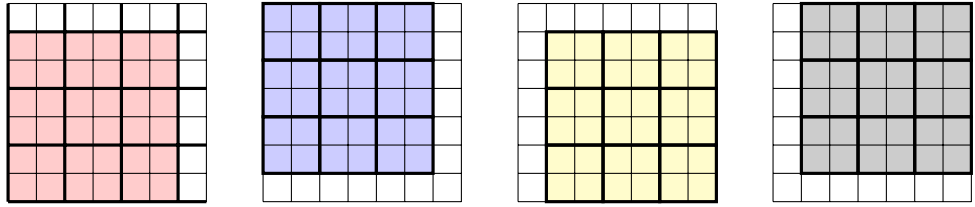}
    \caption{Non-overlapping coloring for vertex patches on regular meshes. In each color we strive to obtain a tiling of the entire domain. Subsequent colors are obtained by shifting the first patch by one cell in each coordinate direction.}
    \label{fig:patch_coloring}
\Description[Fully described in the text.]{}
\end{figure}

With the choices of $V_j$ and $\overline V_j$ of the clamped kernel, the updates generated by one patch only influence the updates of another patch when they overlap. We can leverage this property to partition the set of patches. In fact, all non-overlapping patches depicted in Figure~\ref{fig:patch_coloring} can be processed simultaneously without causing any conflicts. To summarize this, we introduce the colorized vertex-patch smoother in Algorithm~\ref{alg:Multiplicative_Schwarz}. The impact of this localization on the performance of the vertex-patch smoother with continuous elements on the CPU and GPU is studied in~\cite{Wichrowski2023} and~\cite{Cui2023}, respectively. 

\begin{algorithm}
\caption{Colorized Vertex-Patch Smoother $S_\ell(x_\ell,b_\ell)$.}
\label{alg:Multiplicative_Schwarz}
\begin{algorithmic}[1]
\For{$c=1,\dots,n_{c}$} \Comment{Sequential loop over colors}
\For{$j \in\mathcal{J}_c$} \Comment{parallel loop inside color}
\State $r_j \gets R_j b - \overline A_j \overline R_j x_j$ \Comment{local residual}
\State $x_\ell \gets x_\ell + R^T_j A_j^{-1} r_j$ \Comment{local solver}
\EndFor
\EndFor
\end{algorithmic}
\end{algorithm}

Since the domain of dependence of the full kernel extends further to the neighbors of the patches, the structure is more complicated in this case.
In addressing these issues, we decided for updating residuals globally.

\subsection{Fast Diagonalization of Local Solvers}
\label{sec:fd}
In the smoothing operator discussed in the previous section, the direct computation of the local inverse $A^{-1}_j$ requires $\mathcal{O}(k^{3d})$ arithmetic operations,
whereas the subsequent matrix-vector multiplication involves $\mathcal{O}(k^{2d})$ operations.
This stands in stark contrast to the more efficient matrix-free evaluation $A_j v_j$, which demands only $\mathcal{O}(dk^{d+1})$ operations.
This situation is aggravated by the fact that $k$ for the vertex patch is almost twice as large as for a single cell.
Consequently, in this paper, we consider the fast diagonalization method~\cite{lynch1964direct} as an approach to efficiently compute the inverse of $A_j$
\begin{equation}
\label{eq:inverse2d}
    A_j^{-1} = S_1 \otimes S_0 (\Lambda_1 \otimes I + I \otimes \Lambda_0)^{-1} S_1^T \otimes S_0^T
\end{equation}
in two dimensions and
\begin{equation}
\label{eq:inverse3d}
    A_j^{-1} = S_2 \otimes S_1 \otimes S_0 (\Lambda_2 \otimes I \otimes I + I \otimes \Lambda_1 \otimes I + I \otimes I \otimes \Lambda_0)^{-1} S_2^T \otimes S_1^T \otimes S_0^T
\end{equation}
in three dimensions, where $S_d$ is the matrix of eigenvectors to the generalized eigenvalue problem in the given tensor direction $d$:
\begin{equation}
\label{eq:fast_inverse}
    L_dS_d=M_dS_d\Lambda_d , \quad d=0,...,\mathrm{dim},
\end{equation}
and $\Lambda_d$ is the diagonal matrix representing the generalized eigenvalues $\lambda_i$. 

\subsection{Efficiency and Convergence Results}
We have introduced local solvers with three different patch spaces $V_j$, one of them with inconsistent residuals. Hence, we investigate the influence of the local subspace on the convergence of the overall multigrid method.
\begin{table}[tp]
\centering
\caption{Iteration count $\nu$ for vertex patch smoother with \textit{Full} kernel. GMRES solver with relative accuracy $10^{-8}$ preconditioned by multigrid. Entries “---” are not computed due to excessive memory.}
\begin{tabular}{c|ccccc}
\hline
$L$ &$\mathbb{Q}_{3}$ & $\mathbb{Q}_{4}$ & $\mathbb{Q}_{5}$ & $\mathbb{Q}_{6}$ & $\mathbb{Q}_{7}$ \\
\hline
2 & 3.4 & 2.9 & 2.8 & 2.6 & 2.4 \\
3 & 3.7 & 3.2 & 2.8 & 2.6 & 2.3 \\
4 & 3.6 & 3.1 & 2.7 & 2.5 & 2.2 \\
5 & 3.4 & 2.8 & 2.5 & 2.4 & 2.0 \\
6 & 3.0 & 2.6 & 2.4 & 2.2 & 1.8 \\
7 & 2.8 & 2.5 & 1.9 & 1.9 & --- \\
\hline
\end{tabular}
\label{tab:exact-3d}
\end{table}

\begin{table}[tp]
\centering
\caption{Iteration count $\nu$ for vertex patch smoother with \textit{Dirichlet} and \textit{clamped} kernels. GMRES solver with relative accuracy $10^{-8}$ preconditioned by multigrid. Entries “---” are not computed due to excessive memory.}
\begin{tabular}{c|ccccc|ccccc}
&\multicolumn{5}{c|}{Dirichlet kernel}
&\multicolumn{5}{c}{clamped kernel}
\\
\hline
$L$ 
&$\mathbb{Q}_{3}$ & $\mathbb{Q}_{4}$ & $\mathbb{Q}_{5}$ & $\mathbb{Q}_{6}$ & $\mathbb{Q}_{7}$
&$\mathbb{Q}_{3}$ & $\mathbb{Q}_{4}$ & $\mathbb{Q}_{5}$ & $\mathbb{Q}_{6}$ & $\mathbb{Q}_{7}$ \\
\hline
2 & 4.9 & 4.4 & 4.8 & 4.8 & 5.0
  & 18.6 & 11.1 & 7.9 & 6.6 & 5.4\\
3 & 5.7 & 5.9 & 6.1 & 6.2 & 6.5
  & 22.7 & 12.2 & 7.7 & 6.6 & 5.2\\
4 & 5.6 & 6.1 & 6.1 & 6.1 & 6.4
  & 21.7 & 11.6 & 7.4 & 5.8 & 4.8\\
5 & 5.4 & 5.7 & 5.8 & 5.8 & 6.2
  & 19.0 & 10.7 & 6.9 & 5.7 & 4.7\\
6 & 5.3 & 5.5 & 5.6 & 5.6 & 5.9
  & 17.6 & 9.7 & 6.6 & 5.5 & 4.6\\
7 & 4.9 & 5.3 & 5.4 & --- & ---
  & 16.4 & 8.8 & 6.0 & --- & ---\\
\hline
\end{tabular}
\label{tab:cf-3d}
\end{table}
Convergence results are presented in Tables~\ref{tab:exact-3d} and~\ref{tab:cf-3d}.
We solve the linear system~\eqref{eq:mat_form} in three dimensions for right hand side $f \equiv 1$ and for varying finest level $L$ to a relative accuracy of $10^{-8}$. We measure efficiency of the preconditioners by reduction of the Euclidean norm of the residual $\|r_n\|$ after $n$ steps compared to the initial residual norm. Since the (integer) number of iteration steps depends strongly on the chosen stopping criterion,
we define the \emph{fractional iteration count} $\nu$ by
\begin{equation*}
    \nu = -8 \log_{10} \bar{r}, \quad \bar{r} = \left(\frac{\|r_n\|}{\|r_0\|} \right)^{1/n}.
\end{equation*}

The results in Table~\ref{tab:exact-3d} indicate that the number of iterations required for convergence is generally independent of the mesh level. Notably, for degrees 4 and higher, convergence is achieved in at most three steps, making the method nearly as efficient as a direct solver.
A comparison between the data in Tables~\ref{tab:exact-3d} and~\ref{tab:cf-3d}, respectively, reveals that both newly introduced local solvers need more iterations for convergence and are thus numerically less efficient. This is in part due to the smaller subspaces $V_j$, which shows clearly in the high iteration counts for the clamped kernel with small polynomial degrees. For $\mathbb Q_3$, the dimension of $V_j$ is reduced from $8^3$ for the full kernel to $6^3$ for the Dirichlet kernel and $4^3$ for the clamped kernel. Hence, the difference is more pronounced for the clamped kernel. Looking at higher order polynomials, the difference is much less pronounced for the clamped kernel. In the case of the Dirichlet kernel, we see that the inconsistency of the residual is a limiting factor.

\begin{figure}[tp]
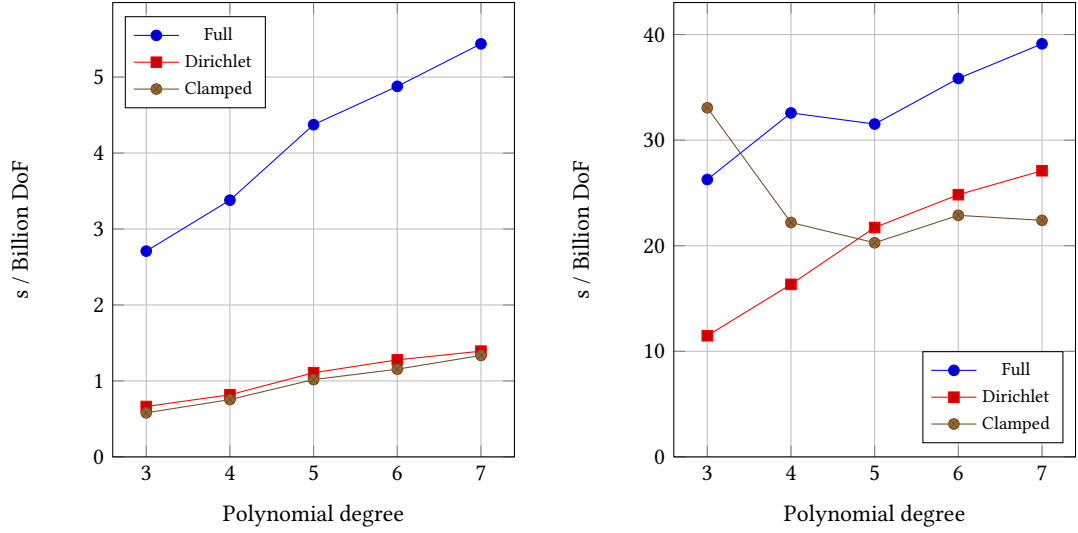

\centering
\includegraphics[width=.45\textwidth]{Figures/Smoothing_poisson.tikz}
\hspace{5mm}
\includegraphics[width=.45\textwidth]{Figures/time_solution_poisson.tikz}
\caption{Comparison of performance of one smoothing step (left) and of the GMRES solver (right) for different kernels with 134--453 million DoFs in 3D.}
\label{fig:Smoothing}
\Description[Fully described in the text.]{}
\end{figure}

However, the benefits of smaller patch sizes, such as enhanced parallelism and data reuse in faster on-chip memory, result in a performance improvement of 4--5 times for a single smoothing step compared to the full kernel, as shown in Figure~\ref{fig:Smoothing} on the left. On the right of the same figure, we show the total time to solution of the preconditioned GMRES method, which clearly shows that the improved memory access of the Dirichlet and clamped kernels outweighs their weakness with respect to convergence rates by far.

\section{GPU Implementation (Single GPU)}\label{sec:Imp_s}

In this section, we explore the GPU implementation of the multigrid method outlined in Section~\ref{sec:GMG}. It is based on the Compute Unified Device Architecture (CUDA) programming model provided by NVIDIA for developing general purpose computing kernels that run on their massively parallel GPUs. 

All numerical experiments in this section are performed on NVIDIA Ampere A100 SXM4 GPU with 80GB of high-speed HBM2e memory for VRAM, which provides up to 2TB/s peak memory bandwidth, The GPU is hosted on a system with two AMD EPYC 7282 16-Core processors.  

When designing the mapping of numerical methods to algorithms executed on GPUs, efficient data access and optimal utilization of the available resources are pivotal for achieving high performance on GPUs. This involves two key considerations: a) Data Layout, which encompasses decisions regarding data requirements, storage location, and data organization or layout; b) Compute Layout, which focuses on identifying the potential parallelism within algorithms and effectively mapping it to the GPU architecture.

\subsection{Data Layout} 

In the matrix-free implementation, the need to store global matrices is eliminated. The only global data storage requirements are for the coefficient vector $x$ and the right-hand side $b$ in~\eqref{eq:mat_form}.
In CUDA programming, optimizing the coalescence of global memory accesses is of utmost significance~\cite{cudaProg}. Employing a lexicographical numbering scheme~\cite{Cui2023} enables an efficient index storage scheme and also ensures that global memory accesses are coalesced whenever possible.
\begin{figure}[tp]
\centering
\includegraphics[width=.65\textwidth]{Figures/dof_layout.tikz}
\caption{Visualization of degrees of freedom layout for $\mathbb{Q}_2$ element in 2D. Left: Global lexicographical numbering. Right: Cell-wise lexicographical numbering.}
\label{fig:dof_layout}
\Description[Fully described in the text.]{}
\end{figure}

However, as illustrated in the left panel of Figure~\ref{fig:dof_layout}, this numbering scheme exhibits two drawbacks. This format constrains problem-solving to a single brick, thereby restricting the choice of computational domains. Moreover, it impedes the potential for MPI parallelism through automated partitioning tools. Achieving global lexical numbering becomes impractical, thus hampering the efficient distribution of computational tasks across multiple processing units.

\begin{figure}[tp]
\centering
\includegraphics[width=.55\textwidth]{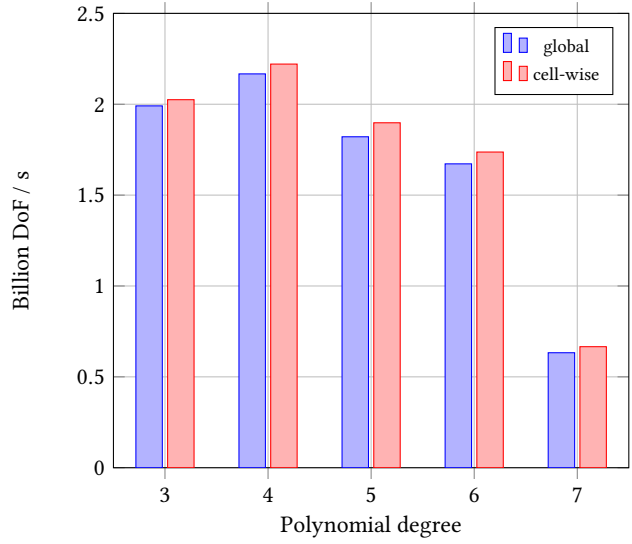}
\caption{Comparison of throughput of different numbering scheme for matrix-free evaluation $Ax$ in 3D. Two different numbering scheme are shown in Figure~\ref{fig:dof_layout}, \textit{global} lexicographical and \textit{cell-wise} lexicographical numbering, respectively.}
\label{fig:l_vs_h}
\Description[Fully described in the text.]{}
\end{figure}

To address these issues, we adopt a cell-wise lexicographical numbering scheme, as demonstrated in the right panel of Figure~\ref{fig:dof_layout}. Given our focus on DG methods, this ordering not only offers improved intuitiveness but also effectively addresses the problems previously mentioned. However, it introduces a new challenge: tensor operations, for example in~\eqref{eq:inverse2d}, necessitate lexicographical ordering of degrees of freedom within a patch.
To overcome this challenge, we implement a renumbering process. After ensuring that the cells within a patch are correctly ordered\footnote{The ordering of cells is established based on the position of the shared vertex within each patch. For instance, in a two-dimensional case, if the shared vertex corresponds to the top-right vertex in a cell, that cell is designated as the first. If the shared vertex corresponds to the top-left vertex in a cell, then the cell is considered the second, with the bottom-right vertex corresponding to the third, and so on.}, we establish a mapping between the two ordering methods for all patches. This fixed mapping is then stored in the constant memory. To validate its suitability, we compare the performance of both ordering methods for computing $Ax$ in Figure~\ref{fig:l_vs_h}. The results demonstrate that, for both lower and higher order, the cell-wise lexicographical numbering scheme exhibits slightly superior performance due to its increased data localization. Consequently, we adopt this numbering format for all subsequent discussions

Now, let us shift our focus to local data on a single patch. The initial step is to establish a global-to-local mapping, as defined by~\eqref{eq:local-solver}, enabling us to extract coefficients associated with the patch from global vectors. With the numbering scheme of Figure~\ref{fig:dof_layout} on the right, it is evident that we can reconstruct the corresponding degrees of freedom of the patch by considering the first degree of freedom from each cell. Maintaining patch indices across different colors demands roughly eight times more memory than storing the coefficient vector $\mathbf{x}$ in three-dimensional situations. By storing only eight degrees of freedom per patch, less than one-eighth of the original memory capacity is needed, especially for higher-order elements, resulting in significant memory savings and enabling the solution of larger-scale problems even with limited resources.

Leveraging the tensor-based basis functions, we can efficiently store the one-dimensional stiffness and mass matrices in shared memory. This reduces the need to access the higher-latency main memory during computations. In our experiments, we found that storing one-dimensional matrices in constant memory greatly increases the need for registers, mainly because the compiler performs various optimizations while compiling the kernel code. It may decide to use registers to store frequently accessed constant memory values to speed up accesses. This, in turn, limits the number of thread blocks that can run simultaneously on the Streaming Multiprocessors (SMs), ultimately impacting overall efficiency.

\subsection{Compute Layout}

The challenge we seek to address revolves around the optimal distribution of tasks for parallel execution on GPUs. In the CUDA programming model, there are two layers of parallelism: parallelism over the thread-blocks and parallelism within each thread-block. It is important to note that communication between blocks is not feasible, while all threads have access to the shared memory within the thread-block. Our design approach is guided by this distinction. We refrain from conceptualizing the entire V-cycle as a single, massive kernel. Instead, we decompose the V-cycle into distinct operations such as $Ax$, the smoothing operation, and the grid transfer operation on individual levels, treating them as individual kernels. These kernels are executed sequentially, with their initiation orchestrated from the host.

As outlined in our algorithm in Section~\ref{sec:mf_op} and Section~\ref{sec:GMG}, the operations within the V-cycle are organized into patches, with each thread block assigned one or more patches. To handle conflicting writes during global vector updates, we make use of atomic operations~\cite{KronbichlerLjunqkvist19}. In pursuit of enhanced parallel execution for smoothing operators, we implement checkerboard coloring, optimizing the algorithm for parallel performance.

Within each thread block, the updates of individual degrees of freedom do not interfere with each other, enabling effective parallelization using a one-thread-per-output strategy. However, due to the limitation of creating a maximum of 1024 threads in a thread block, three-dimensional computations require a single thread to handle updates for multiple (one-column) degrees of freedom. This situation proves advantageous for the scheduler, as it widens the scope of optimization possibilities by involving each thread in a more substantial computational workload, as detailed in~\cite{Cui2023}.

For grid transfer operations, the result is computed cell by cell by looping over the parent-level cells and performing a series of tensor contractions on the input values~\cite{Cui2023}. In the process of prolongation, we first interpolate the coarse cell node values along the first dimension, generating intermediate values. Subsequently, we interpolate along the second dimension to obtain the fine-level node values.
Each thread, corresponding to a coarse-level DoF, is responsible for computing $2^d$ DoF throughout the computation. Consequently, the multigrid transfer can leverage the same loop infrastructure over cells as the matrix-vector products, utilizing the same atomic operations to avoid race conditions.

\section{Compute Optimizations}\label{sec:opt}

\begin{figure}[tp]
\centering
\includegraphics[width=.9\textwidth]{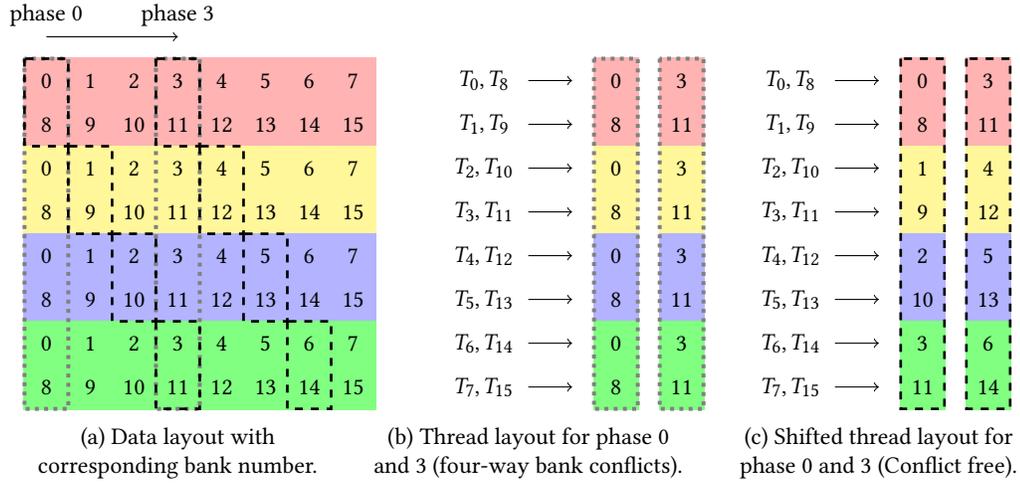}
\caption{Shared memory access pattern for matrix-matrix multiplications with $\mathbb{Q}_3$ element. The dotted and dashed boxes represent the data accessed before and after optimization in phases 0 and 3, along with its arrangement in shared memory (a). When processing the data, threads work from left to right and wrap around. (b) and (c) illustrate the mapping of threads to bank numbers.}
\label{fig:bank_comp}
\Description[Fully described in the text.]{}
\end{figure}

Since we employ basis functions constructed from tensor products of one-dimensional objects, all matrix-vector multiplications can be mapped to matrix-matrix multiplications. These matrix operations have dimensions of $k \times k$ (representing 1D stiffness and mass matrices) and $k \times k^{d - 1}$ (representing input vectors $x_j$ and intermediate results).

\subsection{Bank Conflicts}

\begin{figure}[tp]
\centering
\includegraphics[width=.55\textwidth]{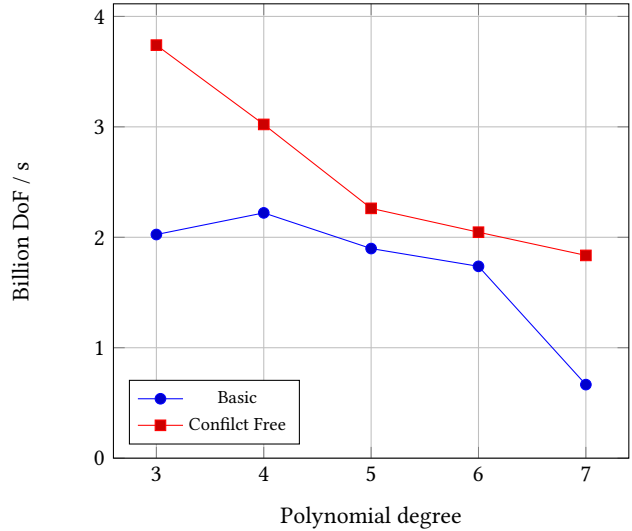}
\caption{Comparison of throughput for shared memory access patterns with and without bank conflicts. Matrix-free evaluation of $Ax$ in 3D.}
\label{fig:Ax_bank}
\Description[Fully described in the text.]{}
\end{figure}

\begin{table}[tp]
    \centering
    \caption{Profiled statistics of bank conflicts per thread block for the different strategies in 3D. The number of bank conflicts are collected based on ``L1 Wavefronts Shared Excessive instruction-level" metric.}
    \label{tab:ncu_bc}
\begin{tabular}{c|ccccc}
\hline
 &$\mathbb{Q}_{3}$ & $\mathbb{Q}_{4}$ & $\mathbb{Q}_{5}$ & $\mathbb{Q}_{6}$ & $\mathbb{Q}_{7}$ \\
\hline
Basic & 1824 & 1040 & 3276 & 2744 & 64974\\
Conflict Free & 0 & 0 & 0 & 0 & 0\\
\hline
    \end{tabular}
\end{table}

In our kernel, the primary computational task is matrix-matrix multiplication. As discussed in our previous paper~\cite{Cui2023}, optimizing the order of operations in matrix multiplication is essential to eliminate bank conflicts as it influences the way we load or store matrices in shared memory. As explained in Figure~\ref{fig:bank_comp}(a) and (b), we access the data in shared memory in a row-major way.
Additional potential for bank conflicts in matrix multiplication arises.
In Figure \ref{fig:bank_comp}(a), we illustrate the bank numbers corresponding to the addresses of shared memory that each thread reads at the start of the computation. Figure~\ref{fig:bank_comp}(b) illustrates the memory and bank numbers associated with different threads in phase 0. Threads 0 and 8 read from bank 0, sharing the same address (broadcast), and this read operation doesn't lead to a bank conflict. The same principle applies to threads 1 and 9. However, threads 2 and 10 also access bank 0, but at a different memory address than threads 0 and 8. Consequently, accessing data in the yellow area occurs after the completion of the red area. This pattern repeats for other threads, resulting in data accesses in phase 0 being split into four sequential completions. By setting the offset based on the number of rows where the threads are located (as shown in Figure~\ref{fig:bank_comp}(c)), we can completely eliminate bank conflicts. 
The impact of bank conflicts on shared memory data loads is examined in~\cite{sun2022dissecting}, revealing that a load instruction with a 4-way conflict incurs a latency of 37ns, which is 1.5 times that of a conflict-free pattern. Figure~\ref{fig:Ax_bank} compares performance in scenarios with and without bank conflicts, clearly demonstrating that avoiding these conflicts markedly improves the algorithm's efficiency. This enhancement is further validated by the \texttt{Nsight Compute}\cite{cudaNcu} results, detailed in Table~\ref{tab:ncu_bc}.

\begin{figure}[tp]
\centering
\includegraphics[width=.55\textwidth]{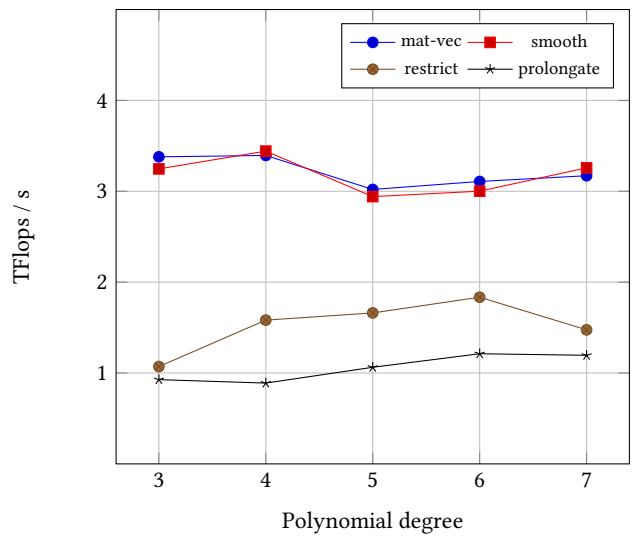}
\caption{Arithmetic throughput in TFlops/s achieved by each components in 3D. The theoretical peak FP64 performance is 8.7 TFlops/s at 1.27 GHz.}
\label{fig:Arithmetic}
\Description[Fully described in the text.]{}
\end{figure}

Recall that we have decomposed the multigrid method into several parts, each performing distinct types of processing, but sharing the same strategy for mapping onto a GPU. Therefore, we assess the performance level achieved by each part of the operator. Figure~\ref{fig:Arithmetic} provides an indication of this performance, based on the number of floating-point operations per second.
Reassuringly, similar flop rates are attained for different orders. However, for the grid transfer operators, a lower performance is expected due to their inclusion of fewer tensor operations and their typical memory-bound nature.
In contrast, the matrix-vector operation and the smoothing operation consistently exhibit high performance in our matrix-free implementation. This is attributed to the more efficient fast diagonalization method employed in the local solver. 
The adoption of the fast diagonalization method approach, however, alters this scenario into matrix-matrix multiplication. This shift enables the algorithm to evolve into a compute bound problem, offering opportunities for achieving higher performance.
The optimal performance for both operation is reached at a polynomial degree of 4, representing 39\% of the theoretical peak.

It is important to note that, firstly, for illustrative purposes, we have opted for double precision; when switching to single precision, the number of banks to 32, necessitating to revisit the permutation to address the bank conflict. However, the underlying principle for resolving this conflict remains consistent. Secondly, to prevent bank conflicts in three dimensions during tensor contraction operations, one must consider not only changes in different rows but also offsets due to variations in the z-direction. In this illustration, we choose the polynomial space $\mathbb{Q}_3$ as an example due to its pronounced characteristics. However, other polynomial degrees may also lead to bank conflicts, as detailed in Table~\ref{tab:ncu_bc}. 
Lastly, it is noteworthy that while padding is commonly employed to mitigate bank conflicts, it introduces a trade-off by increasing pressure on shared memory due to the additional memory usage. Therefore, we do not consider it in this study.

\subsection{Mixed Precision}

The utilization of single-precision arithmetic on GPUs~\cite{komatitsch2010high,goddeke2007performance,gropp2000performance}, either exclusively or partially, provides two major advantages. First, it presents substantial benefits in terms of memory bandwidth. Single-precision data consumes half the memory compared to double precision, effectively doubling the throughput and mitigating memory-related bottlenecks.
The second major factor is the superior peak performance of single precision (FP32) arithmetic compared to double precision (FP64) on Nvidia A100 GPUs. Furthermore, employing single precision results in a smaller memory footprint per thread, lower instruction latency, and consequently, smaller register files and less shared memory usage. This combination of advantages enables the scheduler to allocate more threads to the Streaming Multiprocessors (SM), leading to better occupancy and improved overall performance. In our mixed precision approach, the multigrid V-cycle is fully done in single precision and the outer GMRES iteration is done in double precision. The format is converted when entering and exiting the V-cycle.

\begin{figure}[tp]
\centering
\includegraphics[width=.55\textwidth]{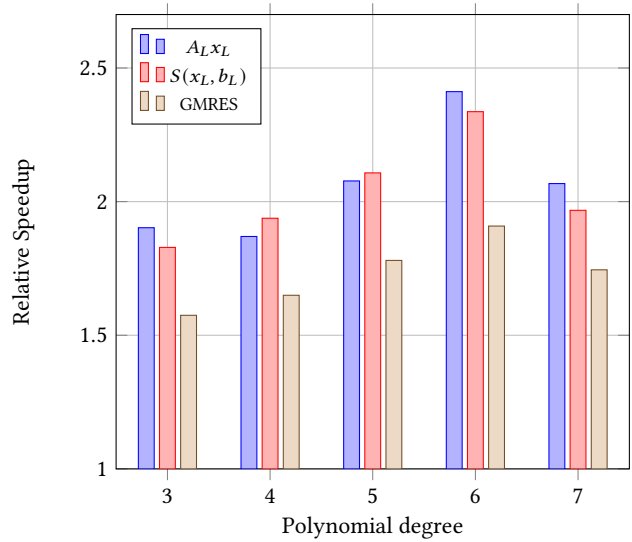}
\caption{Mixed precision performance relative to double precision for different components in 3D. $A_Lx_L$: matrix-vector operation on level $L$; $S(x_L,b_L)$: smoothing operation on level $L$; GMRES: entire GMRES solver.}
\label{fig:speedup_mixed}
\Description[Fully described in the text.]{}
\end{figure}

In Figure~\ref{fig:speedup_mixed}, we evaluate the performance of a GMRES method run in double precision preconditioned by a multigrid V-cycle with one pre- and postsmoothing step in single precision against doing all computations in the multigrid V-cycle in double precision. We find that both the $A_Lx_L$ and $S(x_L,b_L)$ operations not only meet but in some cases exceed a speedup of twice the performance for higher-order elements. This observation is consistent with our theoretical analysis.
However, for low-order elements, the occupancy ratio decreases. This is primarily because we did not modify the algorithm to assign more patches per thread block, but instead ran it in single precision, thus not achieving the expected speedup in this case. Considering that the V-cycle represents the most time-consuming part of the solution process~\cite{Cui2023}, the mixed-precision approach significantly accelerates the solution process, resulting in speedups of over 50\% and, in some cases, even up to 90\%.

\section{GPU Implementation with MPI Parallelization (multiple GPUs)}\label{sec:Imp_m}

The integration of GPU acceleration with distributed processing introduces an additional layer of complexity. This necessitates the preparation of the problem through partitioning and packaging to enable effective solving across multiple GPUs. The partitioning process should meet two key requirements: (i) achieving equally balanced partition sizes and (ii) minimizing ``graph edge cuts" between partitions. To address these requirements, we employ the \texttt{p4est} library~\cite{BursteddeWilcoxGhattas11} provided within the \texttt{deal.II} library~\cite{dealII95}.

\begin{figure}[tp]
\centering
\includegraphics[width=.65\textwidth]{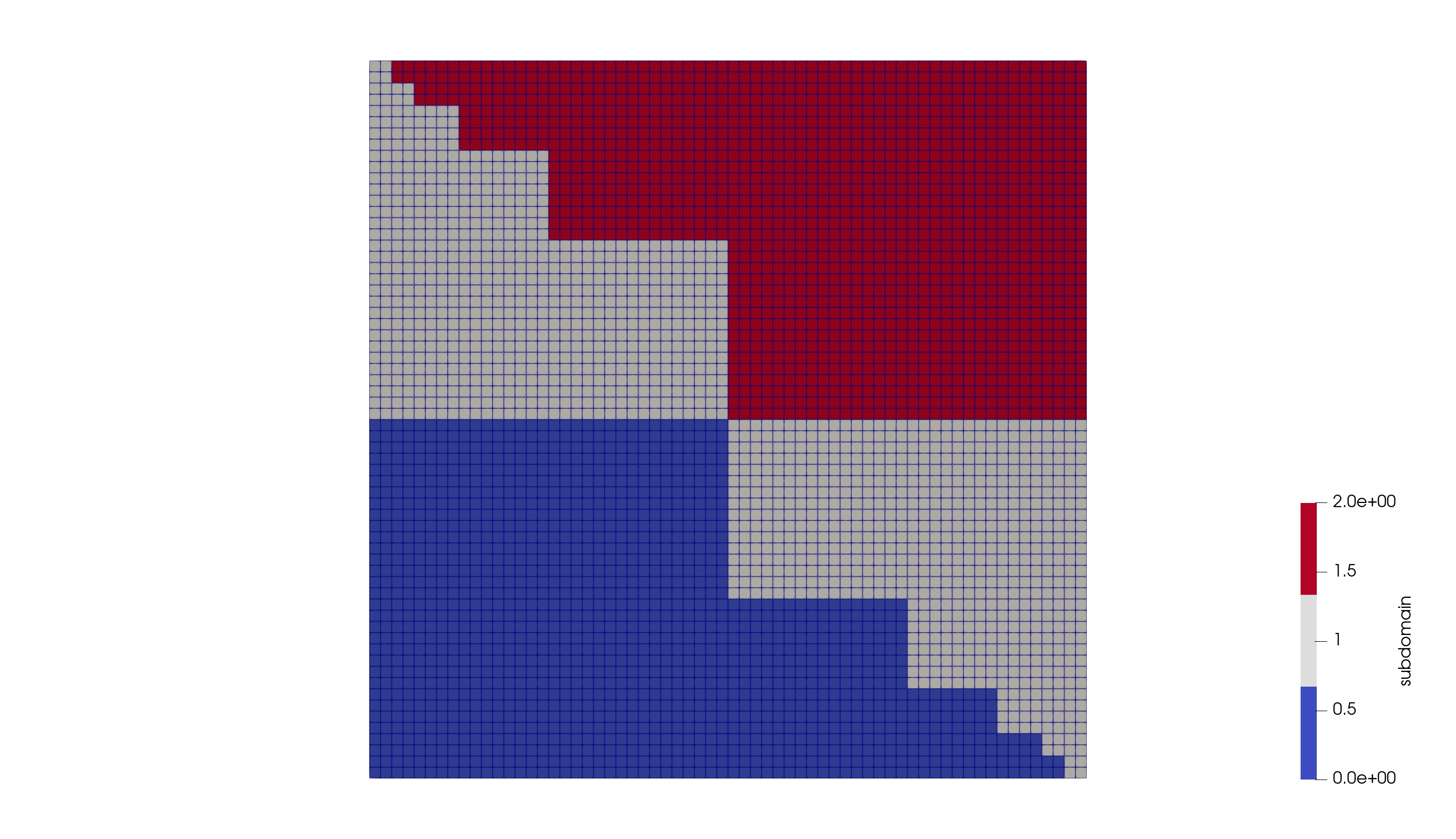}
\caption{A two-dimensional distributed mesh across three processors.}
\label{fig:d_mesh}
\Description[Fully described in the text.]{}
\end{figure}

In the context of finite element computations with MPI parallelism, a common strategy is to decompose the domain, as shown in Figure~\ref{fig:d_mesh}, wherein the cells in the mesh are distributed among the processors. To facilitate the exchange of information between these subdomains, each MPI rank extends its locally owned subdomain by incorporating ghost elements. We use implementations based on \texttt{deal.II}~\cite{BBHK11} and in particular its distributed vector class \texttt{LinearAlgebra::distributed::Vector}. The specific method used for mesh partitioning is not critical, as long as the mesh infrastructure offers the necessary information for uniquely identifying degrees of freedom within the locally owned and ghosted cells.

\subsection{Vertex-patch based Operations}

\begin{figure}[tp]
\centering
\includegraphics[width=.35\textwidth]{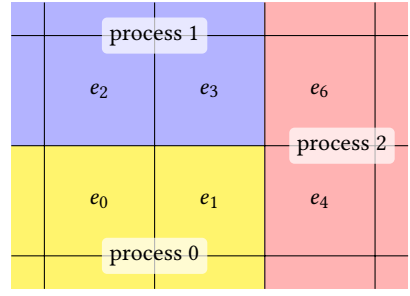}
\caption{Ghost patches with cells owned across by different processes.}
\label{fig:ghost_patch}
\Description[Fully described in the text.]{}
\end{figure}

In an MPI parallelization with the addition of a single layer of ghost cells around locally owned cells, it is possible for cells within a patch to belong to different processors, as depicted in Figure \ref{fig:ghost_patch}. The ownership of such a patch is determined based on two policies: either ensuring that the process that owns the patch includes the fewest ghost elements or assigning ownership to the process with the smallest cell index within the patch.

\begin{figure}[tp]
\centering
\includegraphics[width=.55\textwidth]{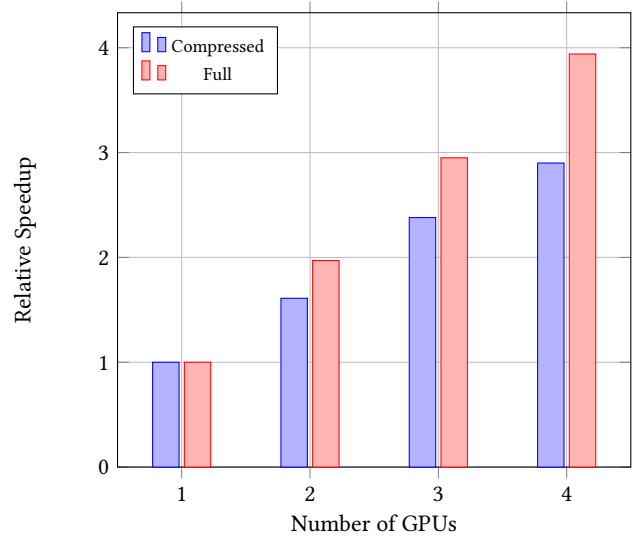}
\caption{Comparison of relative speedup with different storage scheme for ghost patches using $\mathbb{Q}_4$ element in 3D.}
\label{fig:speedup}
\Description[Fully described in the text.]{}
\end{figure}

However, ghost indices in \texttt{deal.II} are stored consecutively at the end of the local index after being sorted according to the global index. This leads to the need for a bisection search when mapping between local and global ghost indices, which can be computationally expensive on the GPU, as each thread may follow a different branch, causing reduced efficiency.
To address this issue, for ghost patches, we explicitly store all the degrees of freedom they contain instead of using the compressed storage method discussed earlier.
By shifting the work of local-global index mapping to the CPU setup stage, we adopt a small device memory overhead as a price for a significant parallel efficiency gain, as demonstrated in Figure~\ref{fig:speedup}.

\begin{figure}[tp]
\centering
\includegraphics[width=.95\textwidth]{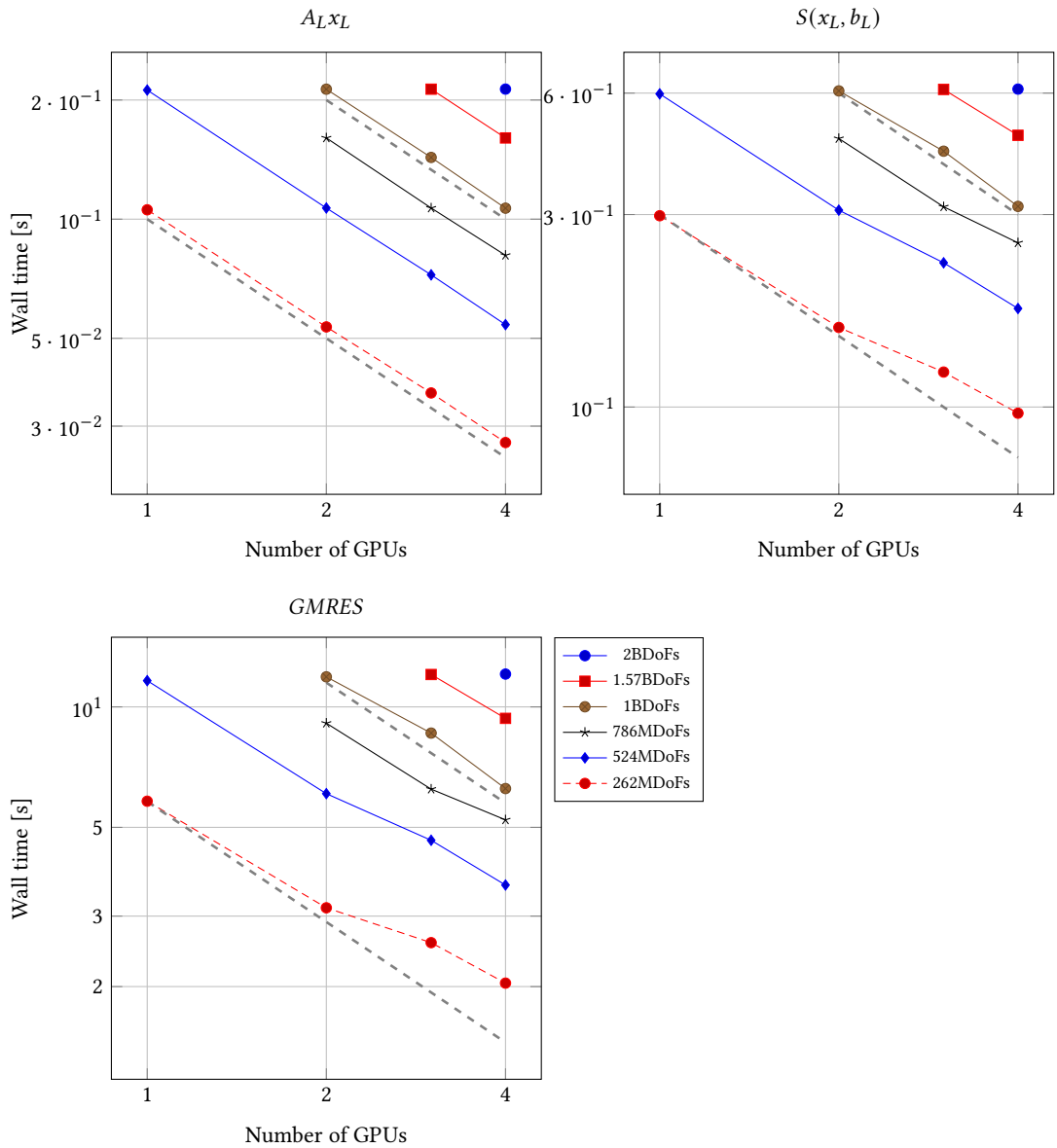}
\caption{Weak and strong scaling analysis of matrix-free evaluation $A_Lx_L$, one smoothing step $S(x_L,b_L)$ and total solution time with $\mathbb{Q}_4$ element in 3D. The greyish-dashed lines serve as references for perfect strong scaling.}
\label{fig:strong_scaling}
\Description[Fully described in the text.]{}
\end{figure}

In our performance and scalability study of the proposed algorithms, we conducted numerical experiments on a server equipped with four A100 GPUs interconnected by NVLink for high-speed multi-GPU communication.

In Figure~\ref{fig:strong_scaling}, we present runtimes for computing the matrix-free operator evaluation $A_Lx_L$, a single smoothing step, and the overall time-to-solution for a polynomial degree of $k = 4$. The greyish-dashed lines serve as references for perfect strong scaling. As expected, the scalability of algorithmic components improves as more degrees of freedom are processed per GPU.
The scaling behavior deteriorates more significantly for smoothing operators, when the problem size is small. This is attributed to the computation of the smoothing operator, which involves updating ghost degrees of freedom. Such updates may introduce ambiguous behavior not allowed in \texttt{deal.II}. Consequently, an additional temporary vector is required, introducing more global communication.  This requirement implies a larger problem size to effectively hide the resulting latency.
\begin{figure}[tp]
\centering
\includegraphics[width=.55\textwidth]{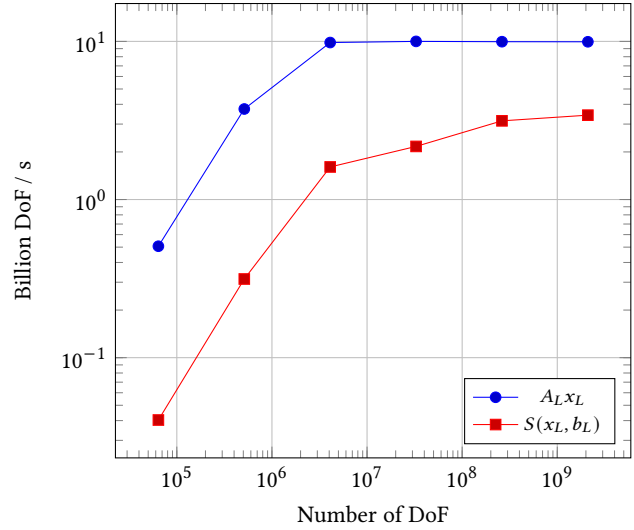}
\caption{Throughput of matrix-free evaluation $A_Lx_L$ and smoothing step $S(x_L,b_L)$ with four GPUs for $\mathbb{Q}_4$ element in 3D.}
\label{fig:throughput_4}
\Description[Fully described in the text.]{}
\end{figure}

In Figure~\ref{fig:throughput_4}, we examine the throughput characteristics of both the matrix-vector product and a single smoothing step concerning problem size, measured in degrees of freedom per second. The data clearly illustrates a rapid attainment of peak performance for the matrix-free evaluation, achieved at a problem size of about 4 million degrees of freedom (DoFs). In contrast, the smoothing operator requires a problem size exceeding 1 billion DoFs to achieve optimal performance. This observation aligns with our prior analysis highlighting the decline in parallel efficiency associated with the smoothing operation.

\begin{figure}[tp]
\centering
\includegraphics[width=.55\textwidth]{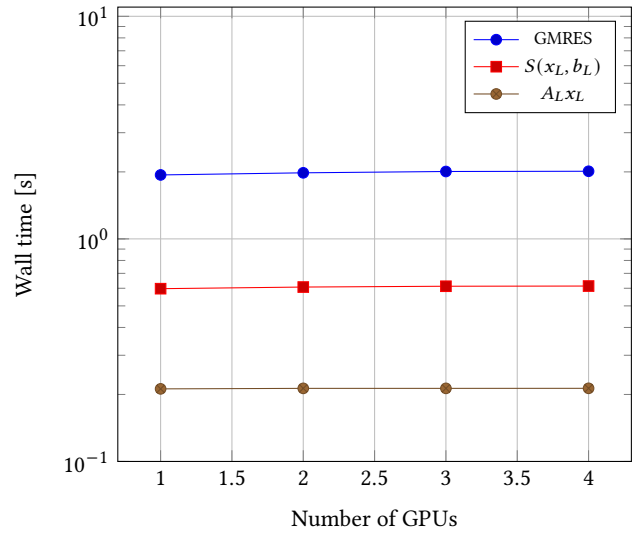}
\caption{Weak scaling results of matrix-free evaluation $A_Lx_L$, one smoothing step $S(x_L,b_L)$ and GMRES solver time per iteration with $\mathbb{Q}_4$ element in 3D.}
\label{fig:weak_scaling}
\Description[Fully described in the text.]{}
\end{figure}

Moving on to weak scalability, we present the single runtimes for $A_Lx_L$, $S(x_L,b_L)$, and GMRES time per iteration in Figure~\ref{fig:weak_scaling}. In weak scaling, we increase the problem size and the number of GPUs, ensuring that the problem size per GPU remains constant. This approach aims to maintain a flat timeline, indicating that the time per operation remains consistent. Finally, detailed convergence results from the experiments are summarized in Table~\ref{tab:parallel_it}. Once again, these results affirm the solver's robustness as refinement increases.

\begin{table}[tp]
\centering
\caption{Iteration count $\nu$ of weak scaling for vertex patch smoother with \textit{Local} kernel using $\mathbb{Q}_3$ element. GMRES solver with relative accuracy $10^{-8}$ preconditioned by multigrid.}
\begin{tabular}{c|cccc}
\hline
 &  \multicolumn{4}{c}{Number of GPUs}\\
\hline
$L$ & 1 & 2 & 3 & 4 \\
\hline
2 & 5.3 & 5.8 & 5.8 & 5.9 \\
3 & 6.2 & 6.3 & 6.3 & 6.3  \\
4 & 6.0 & 6.1 & 6.1 & 6.2 \\
5 & 5.8 & 5.9 & 5.9 & 5.8  \\
6 & 5.6 & 5.7 & 5.9 & 5.9 \\
7 & 5.5 & 5.5 & 5.6 & 5.7 \\
\hline
\end{tabular}
\label{tab:parallel_it}
\end{table}

\section{Conclusion}\label{sec:con}

In this paper, we have described a geometric multigrid method with vertex-patch smoothers for discontinuous Galerkin discretizations on recent Nvidia graphics processors. We studied different choices of local subspaces in order to improve the performance of the smoother.
The methods are specialized for Cartesian meshes and use fast diagonalization for the computation of local inverses. This work has highlighted algorithmic designs to reach high performance and the robustness of preconditioners with respect to the refinement level and polynomial degree was validated on Poisson problems.

We have shown that, with our optimization, matrix-free operations can reach 39 percent of the theoretical peak performance values for the hardware under consideration and demonstrate good scalability with CUDA-aware MPI parallelization. 
To make these speedups accessible, we provide detailed optimization strategies. Notably, the optimization of access patterns to shared memory, eliminating bank conflicts, emerged as the most powerful, resulting in a nearly twofold acceleration.

Furthermore, we have demonstrated the successful integration of a mixed-precision approach within the multigrid method, yielding notable speedups of up to 90\%. A comparative analysis between our GPU implementation and an optimized CPU counterpart running on a cluster revealed that a singular GPU surpasses the performance of an 80-core CPU. This substantiates the pivotal role of localized smoothing operators in achieving heightened computational efficiency.

\begin{acks}
Cu Cui would like to thank China Scholarship Council (CSC NO. 202106380059) for the financial support.
\end{acks}

\bibliographystyle{ACM-Reference-Format}
\bibliography{references}

\appendix

\section{Finite element shape functions and degrees of freedom}

The shape functions of the space $\mathbb{Q}_k$ are defined as tensor products of one-dimensional shape functions, such that we obtain
\begin{gather*}
    \phi_{\alpha\beta}(x,y) = \phi_\alpha(x)\phi_\beta(y),
    \qquad
     \phi_{\alpha\beta\gamma}(x,y,z) = \phi_\alpha(x)\phi_\beta(y)\phi_\gamma(z),
\end{gather*}
in two and three dimensions, respectively. The one-dimensional shape functions $\phi_\nu$ form a basis of the space $\mathbb{P}_k$ of univariate polynomials of degree up to $k$. Hence, the basis of $\mathbb{Q}_k$ is determined after we have chosen the basis of $\mathbb{P}_k$.

For the standard case, we define this basis by Lagrange interpolation on the reference interval $[0,1]$. We choose the set $\{\xi_\mu\}$ of $k+1$ Gauss-Lobatto quadrature points, as these form a basis with a favorable condition number. Then, the basis is uniquely defined by the interpolation conditions
\begin{gather*}
    \phi_\nu(\xi_\mu) = \delta_{\mu\nu}.
\end{gather*}

\end{document}